\documentclass[a4paper,twoside,11pt,reqno]{amsart}

\newcommand{\Ueberschrift}{Twisting by a torsor}

\newcommand{\Kurztitel}{Twisting by a torsor}

\usepackage[english]{babel}

\usepackage[headings]{fullpage}
\usepackage[bottom]{footmisc}

\usepackage{amsmath} 
\usepackage{amssymb} 
\usepackage{amsopn}
\usepackage{amsthm} 
\usepackage{amsfonts} 
\usepackage{mathrsfs}

\usepackage[dvips]{graphicx} \usepackage[usenames]{color}
\usepackage[all]{xy}

\usepackage{dsfont}

\usepackage[active]{srcltx}
\usepackage[pdfpagelabels]{hyperref}

 \hypersetup{
   pdftitle={},
   pdfauthor={},
   pdfsubject={},
  pdfkeywords={},
  colorlinks=true,
  breaklinks=true,
  bookmarksopen=true,
  bookmarksnumbered=true,
  pdfpagemode=UseOutlines,
  plainpages=false}

\newcommand{\bS}{{\mathbb S}}

\newcommand{\dC}{{\mathcal C}}

\newcommand{\dF}{{\mathcal F}}
\newcommand{\dG}{{\mathcal G}}

\newcommand{\dO}{{\mathcal O}}

\newcommand{\dS}{{\mathcal S}}

\DeclareMathOperator{\Isom}{Isom}

\DeclareMathOperator{\Aut}{Aut}

\DeclareMathOperator{\Coh}{Coh}

\DeclareMathOperator{\GL}{GL}

\DeclareMathOperator{\Spec}{Spec}

\DeclareMathOperator{\h}{H}

\DeclareMathOperator{\Gal}{Gal}
\DeclareMathOperator{\Rep}{Rep}

\newtheorem{thm}{Theorem}[section]

\newtheorem{prop}[thm]{Proposition}
\newtheorem{lem}[thm]{Lemma}

\newtheorem{cor}[thm]{Corollary}

\theoremstyle{definition}
\newtheorem{defi}[thm]{Definition}

\theoremstyle{remark}
\newtheorem{rem}[thm]{Remark}

\newenvironment{pro*}[1][Proof]{{\it{#1:}} }{}
\newenvironment{pro**}[1][]{{\it{#1}} }{\hfill $\square$}

\numberwithin{equation}{section}

\author{Michel Emsalem}
\begin{document}

\title[\Kurztitel]{\Ueberschrift}

\date{\today} 

\maketitle
\bibliographystyle{alphaSGA}

\begin{abstract}Twisting by a $G$-torsor an object endowed with an action of a group $G$ is a classical tool. For instance one finds in the paragraph 5.3 of the book \cite{Serre} the description of the "op\'eration de torsion" in a particular context. The aim of this note is to give a formalization of this twisting operation as general as possible in the algebraic geometric framework and to present a few applications. We will focus in particular to the application to the problem of specialization of covers addressed by P. D\`ebes and al. in a series of papers. \\
\end{abstract}


\tableofcontents

\maketitle 
\bibliographystyle{alphaSGA}

\section{Introduction} Twisting by a $G$-torsor an object endowed with an action of $G$ is a classical tool. For instance one finds in the paragraph 5.3 of the book \cite{Serre} the description of the "op\'eration de torsion" in a particular context. We refer to \cite{Breen} for a formulation in a topological framework. The aim of this note is to give a formalization of this twisting operation as general as possible in the algebraic geometric framework and to present a few applications. We will focus in particular to the application to the problem of specialization of covers addressed by P. D\`ebes and al. in a series of papers, which was one of the motivations for writing these notes.

We begin with a section on torsors where we recall classical properties and operations on torsors. It is the opportunity to set the notations. In section 3 we define the twisting operation by a torsor using the cocycle description of a torsor. We state the main properties of this twisting operation in paragraph 4 and recall some classical examples. Paragraph 5 is devoted to the study of the particular case where the twisted objects are themselves torsors. This will lead to the situation we study extensively in the next paragraphs where we consider (ramified) Galois covers that we twist with what we call arithmetic torsors. The aim is to provide a test to know wether there are specializations of the given Galois cover which are isomorphic as torsors to some given arithmetic torsor. 

\section{Torsors} 
\subsection{Left torsors and bitorsors} The aim of this section is to set notations in the description of torsors by cocycles. Let $S$ be a scheme, $G \to S$ a group scheme, $X \to S$ a $S$-scheme and $P \to X$ a left $G$-torsor for some Grothedieck topology. Let $(U_i)_{i \in I}$ be a covering of $X$ trivializing the torsor $P$, and for any $i$ a section $s_i: U_i \to P$. The section $s_i$ induces a trivialization $\theta _i : U_i \times _S G_{|U_i} \to P_{|U_i}$ defined by $\theta _i (x,g)=gs_i(x)$. 

For $i,j \in I$, on the intersection $U_{ij}= U_i \times _S U_j$, one gets  ${s_i}_{|U_{ij}}=g_{ij} {s_j}_{|U_{ij}}$ for some $g_{ij} \in G(U_{ij})$ and the commutative diagram

$$\xymatrix{ P_{|U_{ij}} & U_{ij} \times _S G_{|U_{ij}}\ar[l]_{{\theta _i}_{|U_{ij}}\quad \quad }\ar[d]^\simeq \\
& U_{ij} \times _S G_{|U_{ij}}\ar[lu]^{{\theta _j}_{|U_{ij}}}\\
}$$

where the vertical map is defined by $(x,g) \to (x,gg_{ij})$ which is an isomorphism of trivial torsors. The torsor $P$ is obtained by gluing these trivial left torsors over the $U_{ij}$s.

The $g_{ij}$ clearly satisfy the cocycle condition on $U_{ijk}$ :
$$g_{ij} g_{jk}=g_{ik}.$$

Let us describe the group $H=ad(P)$ of automorphisms of $P$ : locally on each $U_i$ via $\theta _i$, sections of $H$ over $U_i$ are right multiplication by some element 
$g_i\in G(U_i)$. If the $g_i$ are the restriction to $U_i$ of a global section of $H$, they will make the following diagrams commutatve :

$$\xymatrix{ g \ar[r] \ar[d] &gg_i \ar[d] \\
gg_{ij} \ar[r] &gg_{ij}g_j=gg_ig_{ij}}.$$

One gets the relations $$g_i=g_{ij}g_jg_{ij}^{-1}.$$

The $X$-group $H$ is the inner form of $G$ determined by the left $G$-torsor $P$ ; $H$ acts on the right on $P$ by the rule $p.h=h^{-1} (p)$ ($p\in P, h \in H$), and the actions of $G$ and $H$ commute. So $P$ can be viewd as a $(G,H)$-bitorsor. 

The link between left action of $G$ and right action of $H$ can be described by isomorphisms $u_i : H_{|U_i} \to G_{|U_i}$ defined in the following way
$$\forall h \in H_{|U_i} \quad s_ih=u_i(h)s_i.$$ One easily checks that

$$\forall h \in H_{|U_{ij}}\quad u_i(h)=g_{ij} u_j(h) g_{ij}^{-1}.$$

\subsection{Contracted product}
Let $P$ be a $(G,H)$-bitorsor on $X$ and $Q$ a $(H,K)$-bitorsor on $X$ respectively described by $(g_{ij}, u_i)$ and $(h_{ij}, v_i)$, one defines the contracted product $P\wedge ^H Q$ as the quotient of $P \times _X Q$ moded by the relation $(yh,z)=(y,hz)$, pour $h \in H$. 

\begin{lem} \label{produit}In terms of cocycles the contracted product $P\wedge ^H Q$ is described by $(g_{ij}u_j(h_{ij}), u_iv_i)$.
\end{lem}

\proof Let $(U_i)_{i \in I}$ be a trivializing covering for the two torsors with sections $s_i : U_i \to P$ and $t_i : U_i \to Q$; one denotes by $\overline {(s_i,t_i)}: U_i \to P\wedge ^H Q$ the corresponding section. On $U_{ij}$, one gets the equality 
$$\overline {(s_i,t_i)}= \overline {(g_{ij}s_j,h_{ij}t_j)}= \overline {(g_{ij}s_jh_{ij},t_j)}= \overline {(g_{ij}u_j (h_{ij})s_j,t_j)}=g_{ij}u_j (h_{ij})\overline {(s_j,t_j)}.$$ The second formula is obvious.\endproof

\subsection{Right torsors}
The same kind of description holds for a right $H$-torsor $P$.  One defines cocycles with values in $H$ by ${s_i}_{|U_{ij}}={s_j}_{|U_{ij}}h_{ij}$, the $h_{ij}$'s satisfying the cocycle condition
$$h_{ik}=h_{jk}h_{ij}.$$
The group $G$ is the group of automorphisms of the right $H$-torsor, a global section of $G$ being given by local sections $h_i \in H(U_i)$ satisfying $$h_i=h_{ij}^{-1} h_j h_{ij}.$$

One has isomorphisms $v_i :G_{|U_i} \to H_{|U_i}$ satisfying $gs_i=s_iv_i(g)$ and over $U_{ij}$, 
$$v_i=h_{ij}^{-1} v_jh_{ij}.$$

If $(g_{ij}, u_i)$ is the description of $P$ as a left $G$-torsor, its description as right $H$-torsor is given by $(h_{ij}, u_j^{-1})$ where $h_{ij}=u_j^{-1} (g_{ij})$.

\subsection{Inverse of a torsor}
Let $P$ be a $(G,H)$-bitorsor. It is endowed with a right action of $G$ and a left action of $H$ defined by the formulas

$$y\star g=g^{-1}y \quad {\rm and} \quad h\star y=yh^{-1}.$$
Through these actions $P$ is a $(H,G)$-bitorseur that one denotes by $P^0$.

With the preceding notation, local sections $s_i : U_i \to P$, and left cocycle $g_{ij}$ with values in $G$ (resp. right cocycle $h_{ij}$ with values in $H$ ) describing $P$ as a left $G$-torsor (resp. as a right $H$-torsor) defined by the formulas 
$$s_i=g_{ij}s_j= s_j \star g_{ij}^{-1} \quad {\rm and} \quad s_i=s_j h_{ij}=h_{ij}^{-1} \star s_j.$$
Let us denote as before $u_i : H_{|U_i}\to G_{|U_i}$ the isomorphism given by $s_ih=u_i(h)s_i$. The description of $P^0$ in terms of cocycles is given by $(g_{ij}^0,u_i^0)$, where $g_{ij}^0=h_{ij}^{-1}=u_j^{-1} (g_{ij}^{-1})$ and $u_i^0 (g)=u_i^{-1} (g)$. 

Lemma \ref{produit} shows that $P^0$ is an inverse of $P$ for the contracted product.

It is easy to check the following property.

\begin{lem}\label{inverse-du-produit} With the notation of Lemma \ref{produit}, $(P\wedge ^H Q)^0 \simeq Q^0 \wedge ^H P^0$.
\end{lem}
\section{Torsion by a torsor}
\medskip
Suppose we are given a stack ${\dC} \to {\rm Schemas}/S$ with ${\rm Schemas}/S$ endowed with $fpqc$-topology, an affine $S$-group-scheme $G$ locally of finite type anf faithfully flat. Consider the cat\'egory ${\dC}Ê_G (X)$ with objects $\xi $ of ${\dC} (X)$ endowed with a morphism of sheaves $\varphi : G_X \to \underline {\Aut }(\xi )$ ; a morphism from $(\xi , \varphi )$ to $(\xi ' , \varphi ')$ is a morphism $f: \xi \to \xi '$ compatible with $\varphi , \varphi '$, which means that for any morphism of $S$-schemes $\lambda :T \to X$ and for any $g \in G(T)$ the following diagram is commutative:

$$\xymatrix {
\lambda ^\star \xi \ar[d] _ {\varphi (g)}\ar[r] ^{\lambda ^\star f } &\lambda ^\star \xi ' \ar[d] ^{\varphi ' (g) }\\
\lambda ^\star \xi \ar[r]^{\lambda ^\star f }&\lambda ^\star \xi '\\
}$$

\medskip
\begin{thm}\label{thgeneral}
\enumerate{
\item To any right torsor $P \to X$ under $G_X$ is associated a fonctor 
$$\Phi ^P : {\dC}Ê_G(X) \to {\dC} (X)$$

and for any object $\xi$ of ${\dC}Ê_G(X)$ an isomorphism of $\underline {\Aut} ( \xi )$-torsors 

$$\underline{\Isom} _{{\dC} (X)} (\xi , \Phi ^P\xi ) \to P \wedge ^{G_X} \underline {\Aut} ( \xi ) $$

where $ P \wedge ^{G_X} \underline {\Aut} ( \xi )$ is the contracted product of $P$ with $\underline {\Aut} (\xi )$ over $G_X$. %We will denote $\Phi ^P(\xi ) = \xi ^P$.

\item The torsor $P$ determines an inner form $H$ de $G_X$ making $P$ a $(H,G)$-bitorsor and $\Phi ^P$ induces an equivalence of categories ${\dC}_G (X) \simeq {\dC}_{H} (X)$.

\item If $P$ is a $(H,G)$-bitorsor and $Q$ a $(K,H)$-bitorsor, $\Phi ^Q \circ \Phi ^P\simeq \Phi ^{Q\wedge ^HP}$.

\item Conversely if $\xi$ and $\xi '$ are locally isomorphic objects of ${\dC} (X)$, then $P=\underline{\Isom} ( \xi , \xi ')$ is a right torsor under the $X$-group ${\underline {\Aut} (\xi ) }$ and the twisted object $\Phi ^P(\xi )$ of $\xi$ by $P$ is canonically isomorphic to $\xi '$.

}

\end{thm}
\proof
\begin{enumerate}
\item Let $u_i : U_i \to X$, $i \in I$ be a trivializing covering of $X$ for the right torsor $P$ and the corresponding cocycle $g_{ij}\in G(U_{ij})$ and its image $\varphi (g_{ij})= \bar g_{ij} \in \Aut ( \xi _{|U_{ij }})$. The $\bar g_{ij}$'s define descent data for the family of objects $u_i ^\star \xi $ and these descent data are effective. Thus there exists an unique object $\Phi ( \xi ) $ over $X$ endowed with isomorphisms $\theta _i : u_i ^\star \xi \to u_i ^\star \Phi (\xi) $ making all diagrams commutative :

$$\xymatrix{
u_{ij}Ê^\star \xi  \ar [r]^{ { \theta _i} _{|U_{ ij } } } \ar[d] _{\bar g _{ ij } }&  \Phi (\xi ) _ { |U_{ ij } }  \\
u_{ij}Ê^\star \xi \ar [ru]_ { { \theta _j}  _ { |U_ { ij }  }   } & \\
}.$$

One checks that the object $\Phi (\xi )$ doesn't depend on the trivializing covering neither on the chosen representative $g_{ij}$.

\medskip
Let $f : (\xi _1, \varphi _1) \to ( \xi _2, \varphi _2)$ be a morphism between two objects of ${\dC}Ê_G (X)$. For all $i,j$, the following diagrams are commutative :

$$\xymatrix{{u_i ^\ast \xi _1}_{| U_{ij}}  \ar[r]^{f_{|U_{ij}}}  \ar[d]_{\varphi _1 (g_{ij})} &{u_i ^\ast \xi _2}_{| U_{ij}}\ar[d]^{\varphi _2 (g_{ij})}\\
{u_j^\ast \xi _1}_{| U_{ij}}\ar[r]^{f_{|U_{ij}}}  & {u_j^\ast \xi _2}_{| U_{ij}}\\}$$

which means that the $\Phi ( f)=(f_{|U_i})_i$ define a morphism $\Phi ( \xi _1) \to \Phi (\xi _2)$. This defines a functor $\Phi $.

\medskip
Let $\bar h_{ij}$ be another cocycle with values in $\underline {\Aut} (\xi )$ defining another object $\Phi ' (\xi )$. A morphism $\lambda : \Phi '(\xi ) \to \Phi (\xi )$ is the data of morphisms $\lambda _i \in Hom (u_i ^\star \xi , u_i ^\star \xi )$ making the following diagrams commutative:

$$\xymatrix{
u_{ij}^\star \xi \ar[d]^{\bar h_{ ij } } \ar[r] ^ { {\lambda _i }_{ |U_{ ij } } } &u_{ij}^\star \xi  \ar[d]^{\bar g_{ij}} \\
u_{ij}^\star \xi \ar[r] ^ { {\lambda _j }_{ |U_{ ij } } }&u_{ij}^\star \xi \\
}$$

In particular if $h_{ij}$ is the trivial cocycle, $\Phi ' (\xi ) = \xi $ and the preceding diagrams resume to 

$$ {\lambda _j}_{ |U_{ ij } } \circ  {\lambda _i}_{ |U_{ ij } } ^{-1}= \bar g_{ij} $$

which mean that the family $(\lambda _i)$ is a section of the torsor $P \wedge ^{G_X} \underline {\Aut} ( \xi ) $, which corresponds to the image $\bar g_{ij} $ of the cocycle $g_{ij}$ by the morphism $\varphi $. 

This shows a one to one correspondence between section of $P \wedge ^{G_X} \underline {Aut} ( \xi ) $ on $X$ and sections on $X$ of $\underline {\Isom }(\xi , \Phi \xi )$. The same holds over any $T \to X$, which proves the isomorphism of sheaves $\underline{\Isom} _{{\dC} (X)} (\xi , \Phi ^P\xi ) \to P \wedge ^{G_X} \underline {\Aut} ( \xi ) $.
\medskip

\item The inner form $H$ de $G_X$ is obtained by gluing  $u_i ^\ast G$ with $u_j ^\ast G$ over $U_{ij}$ by conjugation by $g_{ij}$. The fact that $H$ acts on $\Phi ( \xi )$ where $\xi $ is an object of ${\dC} _G (X)$ is more or less tautologic: suppos that we are given a section $(g_i)_i$ de $H$, which means sections $g_i \in G(U_i)$ satisfying far all $i,j$
$$(\ast ) \quad {g_j}_{|U_{ij}}\circ g_{ij}= g_{ij} \circ {g_i}_{|U_{ij}}$$

then we have the following commutative diagrams :

$$\xymatrix{{u_i ^\ast \xi }_{| U_{ij}}  \ar[r]^{{g _i}_{|U_{ij}}}  \ar[d]_{\varphi  (g_{ij})} &{u_i ^\ast \xi }_{| U_{ij}}\ar[d]^{\varphi  (g_{ij})}\\
{u_j^\ast \xi }_{| U_{ij}}\ar[r]_{{g _j}_{|U_{ij}}}  & {u_j^\ast \xi }_{| U_{ij}}\\}$$
which associates to the section of  $H$ an automorphism of $\Phi ( \xi )$.

To check that for a morphism $f : \xi _1 \to \xi _2$ in $ {\dC} _G (X)$ its image $ \Phi (f)$ commutes with the action of $H$, one can check it locally on the $U_i$ where it is obvious. 

\medskip

The fact that $\Phi$ is an equivalence is a consequence of point (3)  which implies that $\Phi ^{P^0}$ is a quasi-inverse of $\Phi ^P$.
\medskip
\item Point (3) is a consequence of the commutative diagram 

$$\xymatrix{\xi _{|U_{ij}} \ar[d]_{\bar g_{ij}}\ar[r]^{\theta _i}&\Phi ^P(\xi )_{|U_{ij}}\ar[d]_{\bar h_{ij}}\ar[r]^{\omega _i}&\Phi ^Q(\Phi ^P(\xi ))_{|U_{ij}} \\
\xi _{|U_{ij}} \ar[d]_{\overline{ u_j(h_{ij})}}\ar[ru]^{\theta _j}& \Phi ^P(\xi )_{|U_{ij}}\ar[ru]^{\omega _j}&\\
\xi _{|U_{ij}}\ar[ru]^{\theta _j}&&\\
}$$
and of Lemma \ref{produit} which gives the description by cocycles of the contracted product.
\medskip
\item In the other direction the first assertion is obvious.

Let $U_i$ be a covering of $X$ with isomorphisms $\lambda _i : u_i ^\star \xi \to u_i ^\star \xi '$. The cocycle defining the torsor $P=\underline {Isom} (\xi , \xi ' )$ is $\bar g _{ij} = {\lambda _j }_{|U_{ij}}^{-1} \circ  {\lambda _i }_{|U_{ij}}$. The following diagrams are commutative
$$\xymatrix{
u_{ij}^\star \xi \ar[d]^{\bar g_{ ij } } \ar[r] ^ { {\lambda _i }_{ |U_{ ij } } } &u_{ij}^\star \xi  \ar[d]^{Id} \\
u_{ij}^\star \xi \ar[r] ^ { {\lambda _j }_{ |U_{ ij } } }&u_{ij}^\star \xi \\
}$$

which proves that $\xi '$ is obtained from $\xi $ by the descent data $\bar g _{ij}$, In other words $\xi ' =\Phi ^P( \xi )$.

\end{enumerate} \endproof
\begin{rem} From the proof one clearly gets that if $P$ the trivial torsor in point (2) $H \simeq G$ and the equivalence of category $\Phi ^P$ is the identity of ${\dC} _G(X)$.
\end{rem} 
 \begin{rem} Given a right $G$-torsor $\pi : P \to X$, $\pi ^\ast $ is a functor from $\dC _G (X)$ to the category of $G$-equivariants objects of $\dC (P)$ (see \cite{Vis}, 3.8, for the definition of $G$-equivariants objects). On the other hand, as explained in \cite{Vis}, Theorem 4.46, the $G$-torsor $\pi : P \to X$ defines an equivalence of categories $F_P$ between $G$-equivariant objects of $\dC (P)$ and $\dC (X)$. The functor $\Phi ^P$ of Theorem \ref{thgeneral} is $\Phi ^P=F_P \circ \pi ^\ast$. 
 
 \end{rem}

\section{Properties of twisting by a torsor}
\subsection{Moding by a normal subgroup}
\begin{prop}\label{quotient} Suppose we are given a normal $S$-sub-group $K$ of $G$. Denote by $\Theta : {\dC} _{G/K} (X) \to {\dC} _G(X)$ the obvious functor. For any $G$-torsor $P$ over $X$ and any object $\xi $ of $ {\dC} _{G/K} (X)$, 
$$\Phi ^{P/K}(\xi ) \simeq \Phi ^P ( \Theta ( \xi ) ).$$
\end{prop}

\proof The statement is a consequence of this simple remark: with the notation of the proof of Theorem \ref{thgeneral}, if $g_{ij}$ denotes the cocycle associated to the torsor $P$ and the covering $U_i$, the cocycle associated to the $G/K$-torsor $P/K$ is $\theta (g_{ij})$ where $ \theta : G \to G/K$ is the canonical morphism. \endproof
\subsection{Image by a morphism of stacks}
Let us consider a functor of $S$-stacks $u : {\dC}_1\to {\dC}_2$ over the category of $S$-schemes. Given an object $\xi $ of ${\dC} _{1,G} (X)$ one still denotes by $u$ the morphism $u :\underline { \Aut }_X (\xi ) \to \underline {\Aut }_X (u(\xi ))$. If $ \varphi : G_X \to \underline { \Aut }_X ( \xi )$ is the structural morphism of $\xi$, $\psi = u \circ \varphi : G _X \to \underline { \Aut } _X ( u( \xi ))$ is the structural morphism of  $u( \xi )$ which is an object of ${\dC}_{2,G} (X)$. Let $P \to X$ be a right $G$-torsor over $X$. One easily checks the following statement.

\begin{prop}
$$u ( \Phi ^P(\xi )) \simeq \Phi ^P(u( \xi ) )$$
isomorphic in ${\dC}_{2,H} (X)$ where $H=\underline { \Aut } _G(P)$.
\end{prop}

\medskip

\subsection{Base change}
Using notation of section 1, consider a morphism $f : X'\to X$ of $S$-schemes. Let $P \to X$ be a right $G$-torsor, the pull back $f^\star P \to X'$ is a right $G$-torsor over $X'$.

\begin{prop}\label{changement-de-base} For any object $\xi $ of ${\dC} _G (X)$,

$$\Phi ^{f^\star P}(f^\star ( \xi )) \simeq f^\star (\Phi ^P(\xi ))$$ 
\end{prop}

\proof If we are given effective descent data over $X$, i.e. objects $\xi _i$ on the opens $U_i$ with isomorphisms $\alpha _{ij} : {\xi _i}_{|U_{ij}} \to  {\xi _j}_{|U_{ij}}$ satisfying descent conditions and defining an object $\xi $ on $X$, one can pull them by $f: X' \to X$ and one gets descent data for the objects $f ^\star \xi _i$ on the opens $f^{-1} U_i $ with isomorphisms $f^\star \alpha _{ij} : {f^\star \xi _i}_{|f^{-1} U_{ij}} \to  {f^\star \xi _j}_{|f^{-1} U_{ij}}$ which define the object $f^\ast \xi$. In particular the right $G$-torsor $P$ is defined by gluing together trivial right $G$-torsors $G_{U_i}$ on the opens $U_i$ with descent data given by the cocycle $g _{ij} \in G( U_{ij})$ (which defines an automorphism of the trivial right $G$-torsor $G_{U_{ij}}$ by left multiplication by $g_{ij}$). Pulling these data by $f$ one gets the left multiplication by $g_{ij} \circ f$ over $f^{-1} (U_{ij})$. Thus the torsor $f^\ast P$ is defined by the cocycle $(g_{ij}Ê\circ f)$.

Fix a right $G$-torsor $P$. Let $\xi $ be an object of ${\dC} _G (X)$; we will denote $\overline g$ the image of an element $g \in G$ by the structural morphism $G_X \to \underline{\Aut} (\xi )$. The twisted object $\Phi ^P (\xi )$ is defined by descent data provided by the  $\overline {g_{ij}}$'s acting on the $\xi _{U_{ij}}$'s. On the other hand $f^\star (\Phi ^P (\xi ))$ is defined by descent data $f^\ast (\bar g_{ij})$ on the $f^\star (\xi _{U_i})= f^\star (\xi ) _{f^{-1} (U_i) }$. The fact that $G_X \to \underline{\Aut} (\xi )$ is a morphism of sheaves implies that $f^\ast (\bar g_{ij})= \overline {g_{ij} \circ f}$. This proves that $f^\star (\Phi ^P (\xi ))$ is defined by the local objects $ f^\star (\xi ) _{f^{-1} (U_i) }$ with descent data $ \overline {g_{ij} \circ f}$. This is the definition of $\Phi ^{f^\star P}(f^\star ( \xi ))$.  \endproof

\subsection{Classical examples}
\begin{enumerate}
\item One knows that the isomorphism classes of vector bundles of rank $n$ on a scheme $X$ are parametrized by $\h ^1 (X, \GL _n)$. The pointed set $\h ^1 (X, \GL _n)$ parametrizes also the $\GL _n$-torsors over $X$. The correspondence between vector bundles and $\GL _n$-torsors can be described in terms of twisting by a torsor. In one direction to a vector bundle $\dF$ of rank $n$ corresponds the $\GL _n$-torsor $P=\underline { \Isom } (\dO _X ^{\oplus n}, \dF)$. In the other direction to the $\GL _n$-torsor $P$ corresponds the twisted vector bundle $\Phi ^P (\dO _X^{\oplus n})$.
\item Let $G \to \Spec (k)$  be an affine group scheme over a field $k$ and consider the category $\Rep _{k} (G)$ of representations of $G$ on finite dimensional $k$-vector spaces. It is a neutral tannakian category with the forgetful functor $i$. For any $k$-scheme $\theta : X \to \Spec (k)$ one may consider the functor $\theta ^\ast \circ i: \Rep _{k} (G) \to \Coh (X)$ with values in the category $\Coh (X)$ of coherent sheaves on $X$. There is a one to one correspondence between fiber functors $F: \Rep _{k} (G)\to \Coh (X)$ and $G$-torsors : given a fiber functor $F: \Rep _{k} (G)\to \Coh (X)$ one associates $\underline {\Isom} ^\otimes (\theta \circ i, F)$ which is a right torsor under the automorphism group of the forgetful functor which is $G$; in the other direction to a right $G$-torsor $P \to X$ one associates the twisted by $P$ of the forgetful functor $\Phi ^P (\theta \circ i)$. This correspondence in the case of finite group schemes $G$ underlies the tannakian construction by Nori of the fundamental group scheme \cite{Nor2}. 
\item When $ X = \Spec (k)$ one recovers the equivalence of categories between the gerbe of fiber functors of the category $\Rep _{k} (G)$ and the gerbe $BG_k$ of $G$-torsors over a $k$-scheme.
More generally a gerbe $\dG \to \Spec (k)$ over field $k$ bound by a $k$-group scheme $G$ is neutral if and only if it is equivalent to the category $BG_k$. Indeed if $\dG$ is neutral and $\xi _0$ is a section $\xi _0 : \Spec (k) \to \dG$, to any section $\xi : X \to \dG $ over a $k$-scheme $X$ is associated the $G$-torsor $P = \underline {\Isom}_X (\xi _0, \xi )$. In the other direction to any $G$-torsor $P \to X$ corresponds the twisted section $\Phi ^P ( \xi _0)$.
\item  In the paragraph 5.3 of "Cohomologie Galoisienne" (\cite{Serre}) J.-P. Serre describes the twisting operation in the category of $H$-sets, where $H$ is a profinite group. The groups $G$ and the right $G$-torsors are groups and torsors in the category of $H$-sets. This leads to an interpretation of the cohomology group $\h^1(H,G)$ in terms of $G$-torsors in the category of $H$-sets, and in the situation of a subgroup $G<G_1$ to an interpretation of the fibers of the natural map of pointed sets $\h^1 (H,G) \to \h^1 (H,G_1)$.
\item More generally coming back to the general situation of Theorem \ref{thgeneral}, consider a morphism of $S$-group schemes $u:G \to G_1$. It leads to a map $\tilde u: \h^1 (X,G) \to \h ^1 (X,G_1)$ sending a $G$-torsor $P$ to the contracted product $P \wedge ^{G}G_1$ through $u$. On the other hand, given a right $G$-torsor $P$, $u$ induces a morphism of group schemes $u': G'=\underline {\Aut} _G P \to G_1' =\underline {\Aut} _{G_1} P_1$ where $P_1 = P \wedge ^G G_1$, and thus a morphism $\tilde u': \h^1 (X,G') \to \h ^1 (X,G'_1)$. The following diagram is commutative:
$$\xymatrix{
\h^1 (X,G) \ar[r]^{\tilde u}&\h^1(X,G_1)\\
\h^1 (X,G')\ar[u]^{-\wedge ^{G'}P}\ar[r]_{\tilde u '}&\h^1 (X,G_1')\ar[u]_{-\wedge ^{G'_1}P_1}\\
}$$ where the vertical maps are bijections. 
Thus $\tilde u ^{-1} (P_1)= \tilde u ^{-1}(\tilde u (P))$ is in one to one correspondence by the left vertical map with the kernel of $\tilde u'$, i.e. the set of right $G'$-torsors $Q'$ such that $Q' \wedge ^{G'} G_1 '$ is the trivial $G_1'$-torsor. In other words the right $G$-torsors $Q$ and $P$ have same images by $\tilde u$ if and only if the image of $Q \wedge ^{G'} P^0$ by $\tilde u'$ is the trivial torsor.
\end{enumerate}

\medskip

\begin{section}{Twisting a $G$-torsor}
Consider the situation of Theorem \ref{thgeneral} with $\dC$ being the category of right $G$ torsors and let $R$ be a right $H$-torsor over $X$. The category $\dC _H(X)$ contains the category of $(H,G)$-bitorsors over $X$ and if $P$ is an element of $\dC _H (X)$, one may consider the twisted object $\Phi ^R (P)$ which is a right $G$-torsor over $X$. 

\begin{cor}\label{torsion-de-torseur} In the situation of Theorem \ref{thgeneral}, for any $(H,G)$-bitorsor $P$, $$\Phi ^R (P)\simeq R\wedge ^H P$$ as right $G$-torsors and $$\underline {\Isom }_G(P, \Phi ^R (P)) \simeq R$$ as right $H$-torsors.
\end{cor}

\proof This is a consequence of points (3) and (1) of Theorem  \ref{thgeneral}. \endproof

\begin{cor}\label{isomorphisme} Let $P$ and $Q$ be right $G$-torsors.
$$\underline{\Isom }_G(P,Q) \simeq \Phi ^Q(P^0)$$ as right $\underline{\Aut }_G (P)$-torsors.
\end{cor}
 \proof Apply Corollary \ref{torsion-de-torseur} to $Q=\Phi ^R (P)=R \wedge ^H P$, where $H=\underline{\Aut }_G (P)$, which is equivalent to $R \simeq Q\wedge ^G P^0\simeq \Phi ^Q (P^0)$.
 \endproof
\begin{defi} The object $\Phi ^Q (P^0)$ will be referred as $P$ twisted by $Q$.
\end{defi}

\begin{rem} In view of Lemma \ref{inverse-du-produit} considering $P$ twisted by $Q$ resumes to considering $Q$ twisted by $P$ (they are inverse of each other). 
\end{rem}
\begin{cor} Let $G_d$ be the trivial $G$-torsor. For any right $G$-torsor $P$, $P^0\simeq \underline { \Isom }_G (P,G_d)$. 
\end{cor}

\proof This is an immediate consequence of Corollary \ref{isomorphisme} with $Q=G_d$.
\endproof

Consider again a $(K,G)$-bitorsor $Q$ and a $(L,G)$-bitorsor $P$. The aim of what follows is to give a description is terms of cocycles and descent data of the object $\Phi ^Q (P^0)=Q\wedge ^GP^0$. The right $G$-torsor $Q$ is defined by $(g_{ij},v_i)$ where $g_{ij} \in G_{U_{ij}}$ and $v_i : K_{U_i} \to G_{U_i}$; the right $G$-torsor $P$ is defined by $(g'_{ij},v'_i)$ where $g'_{ij} \in G_{U_{ij}}$ and $v'_i : L_{U_i} \to G_{U_i}$. So the left $G$-torsor $P^0$ is defined by $({g_{ij}'}^{-1}, v'_i)$ and as left $K$-torsor, $Q$ is described by $(v_j^{-1} (g_{ij}),v_i^{-1})$. The formula of the contracted product gives a description of the left $K$-torsor $Q\wedge ^GP^0$ : $(v_j^{-1}(g_{ij}) v_j^{-1} ({g'_{ij}}^{-1}), v_i^{-1} v'_i)$. We state the result :

\begin{prop} Let $Q$ be a $(K,G)$-bitorsor and $P$ a $(L,G)$-bitorsor described as right $G$-torsors by the cocycles $(g_{ij},v_i)$, where $g_{ij} \in G_{U_{ij}}$ and $v_i : K_{U_i} \to G_{U_i}$, and $(g'_{ij},v'_i)$ where $g'_{ij} \in G_{U_{ij}}$ and $v'_i : L_{U_i} \to G_{U_i}$. The left $K$-torsor $\Phi ^{Q} (P^0)$ is described by $$(v_j^{-1}(g_{ij}{g'_{ij}}^{-1}), v_i^{-1} v'_i).$$
\end{prop}

In terms of descent data $P^0$ is obtained by gluing the trivial left $G$-torsors $G_{U_i}$ over $U_{ij}$ by right multiplication by ${g'_{ij}}^{-1}$ : precisely one has isomorphisms of left $G$-torsors $\varphi _i : P^0 _{|U_i}\to G_{U_i}$ making the following diagrams commutative over $U_{ij}$ :

$$\xymatrix{ P^0 _{|U_{ij}} \ar[r]^{\varphi _i}\ar[dr]_{\varphi _j}& G_{U_{ij}} \ar[d]^{({g'_{ij}}^{-1})_d}\\
& G_{U_{ij}}\\
}$$

where $(g)_d$ denotes the right multiplication by $g$. On the other hand by definition of the twisting of a $G$-object by $Q$ one has isomorphisms $\psi _i : \Phi ^Q (P^0)_{|U_i} \to P^0 _{|U_i}$ which get into the commutative triangles on the $U_{ij}$'s :

$$\xymatrix{
\Phi ^Q (P^0)_{|U_{ij}} \ar[r]^{\psi _i} \ar[rd]_{\psi _j}&P^0_{|U_{ij}}\ar[d]_{g_{ij}}\\
&P^0_{|U_{ij}}\\
}$$

where the vertical map is the action of $g_{ij}$ on $P^0_{|U_{ij}}$.

Finally the following diagrams are commutative over $U_{ij}$ :

$$\xymatrix{
\Phi ^Q (P^0)_{|U_{ij}} \ar[r]^{\psi _i} \ar[rd]_{\psi _j}&P^0_{|U_{ij}}\ar[d]_{g_{ij}}\ar[r]^{\varphi _i}\ar[dr]_{\varphi _j}& G_{U_{ij}}  \ar[d]^{({g'_{ij}}^{-1})_d} \\
&P^0_{|U_{ij}}  \ar[dr]_{\varphi _j}& G_{U_{ij}}  \ar[d]^{({g_{ij}})_g}\\
&&G_{U_{ij}}\\}.$$

We have proved the following statement :

\begin{prop} Let $Q$ be a $(K,G)$-bitorsor and $P$ a $(L,G)$-bitorsor described as right $G$-torsors by the cocycles $(g_{ij},v_i)$, where $g_{ij} \in G_{U_{ij}}$ and $v_i : K_{U_i} \to G_{U_i}$, and $(g'_{ij},v'_i)$ where $g'_{ij} \in G_{U_{ij}}$ and $v'_i : L_{U_i} \to G_{U_i}$. The left $K$-torsor $\Phi ^{Q} (P^0)$ is obtained by gluing the $G_{U_i}$'s over the $U_{ij}$'s by the morphisms
$$({g_{ij}})_g \circ ({g'_{ij}}^{-1})_d = ({g'_{ij}}^{-1})_d \circ ({g_{ij}})_g $$
\end{prop}

\begin{rem} As $P^0$ is locally trivial, $\Phi ^Q (P^0)=Q \wedge ^G P^0$ is locally over $X$ isomorphic to $Q$.
\end{rem}

Here is a consequence of Corollary \ref{isomorphisme} :

\begin{cor}\label{condition-isomorphisme}
Let $u : U \to X$ be a morphism of $S$-schemes. The $U$-scheme  $u ^\star (\Phi ^Q(P^0))$ has a section over $U$ if and only if the torsors $u ^\star  P $ and $u ^\star Q$ are isomorphic over $U$.
\end{cor}

\proof This follows immediately from Corollary \ref{isomorphisme} and from the base change formulas. \endproof
\medskip

In the same way one gets the following.

\begin{cor}\label{condition-isomorphisme-bis}
Let $u : U \to X$ be a morphism of $S$-schemes and $Q' \to U$ be a right $G$-torsor. The $U$-scheme  $\Phi ^{Q'} (u^\ast (P^0))$ has a section over $U$ if and only if the right torsors $u ^\ast P $ and $Q'$ are isomorphic over $U$.
\end{cor}

\medskip

The next property deals with the quotient of $G$ by a normal subgroup.

\begin{prop}\label{quotient-partiel} Let $K \hookrightarrow G$ be a normal subgroup do $G$; suppose we are given a right $G$-torsor $P_1 \to X$, a right $G$-torsor $P_2\to X$, a right $G/K$-torsor $R \to X$ inserted in the commutative diagram with equivariant maps

$$\xymatrix{P_1 \ar[rd]^K \ar[rdd]_G&&P_2 \ar[ld]_K\ar[ldd]^G \\
&R\ar[d]^{G/K} &\\
&X&\\ }.$$
Suppose we are given a section $s: X \to \underline {\Isom }_{G/K,X} (R,R)\simeq \Phi ^R_{G/K,X}(R^0)$ (which is a closed immersion). Denoting $s(X) \in  \Isom _{G/K,X} (R,R) $ the corresponding element the following diagram is cartesian:

$$\xymatrix{
\Phi _{K,R}^{s(X)^\ast P_2}(P_1^0) \ar[r] \ar[d] &s^\ast \Phi _{G,X} ^{P_2}(P_1^0) \ar[r] \ar[d] & \Phi _{G,X} ^{P_2}(P_1^0) \ar[d] \\
R \ar[r] & X \ar[r]_s & \Phi _{G/K,X} ^{R}(R^0) \\
}$$
This holds in particular for the section $s$ corresponding to $s(X) = Id _R$. 
\end{prop}

\proof The diagram can be reformulated in
$$\xymatrix{
\underline {\Isom } _{K,R} (P_1,s(X)^\ast P_2) \ar[r] \ar[d] &s^\ast \underline {\Isom } _{G,X} (P_1,P_2) \ar[r] \ar[d] &\underline {\Isom }_{G,X} (P_1,P_2) \ar[d] \\
R \ar[r] & X \ar[r]_s & \underline {\Isom } _{G/K,X} (R,R) \\
}$$

The right square being cartesian one has to check that the large rectangle is cartesian. The question is local and we may suppose that the torsors are trivial. But $\Isom _{K,R} (P_1,s(X)^\ast P_2)\simeq \Isom _{s(X),R} (P_1, P_2)$ where the last term refers to the set of $X$-isomorphisms from $P_1$ to $P_2$ inducing $s(X)$ on $R$.  \endproof
\end{section}

\medskip

Given two right $G$-torsors $P$ and $Q$ over $X$ Proposition \ref{isomorphisme} gives a tool allowing to test for $u : UÊ\to X$ whether $u^\ast P \simeq u^\ast Q$. In the following construction one answers the question to know, given two morphisms $u, v : U \to X$, whether $u^\ast P \simeq v^\ast Q$. 

Consider the right $G$-torsors $ P_1= P \times _S X \to X \times _S X$ and $Q_1 = X \times _S Q \to X \times _S X$ and the twisted object $Q_1 \wedge ^G P_1^0$.

\begin{prop} Let $(u,v) : U \to X \times _S X$ ; the $G$-torsors $u ^\ast P$ and $v^\ast Q$ are isomorphic if and only if the fibre at $(u,v)$ of $Q_1 \wedge ^G P_1^0$ has a section over $U$. 

Moreover if $v : S \to X$ is a section of the structural morphism $X \to S$, $(1_X \times v)^\ast (Q_1\wedge ^G P_1 ^0) \simeq (v^\ast Q \times _S X)\wedge ^G P^0$ ($P$ twisted by the fibre at $v$ of $Q$). A similar statement holds with $(v \times 1_X)$. 
\end{prop}

\proof This comes from the cartesian diagrams

$$\xymatrix{u^\ast P \ar[d] \ar[r] & P \times _S X \ar[d] \\
U \ar[r]^{u \times v}&X \times _S X\\} \quad {\rm and } \quad  \xymatrix{v^\ast Q \ar[d] \ar[r] & X \times _S Q \ar[d] \\
U \ar[r]^{u \times v}&X \times _S X\\}$$ and from commutation of base change with twisting operation. The proof of second part of the statement uses the same tools.

\endproof

In the particular case where $P=Q$ one recovers the self-twisted cover introduced in \cite{Debes}.

\begin{cor} \label{self-twisted} With the above notation $u ^\ast P$ and $v^\ast P$ are isomorphic if and only if the fibre at $(u,v)$ of the self-twisted torsor $P_1 \wedge ^G P_1^0$ has a section over $U$.
\end{cor}

\section{An example of computation: the case of a Galois extension of fields}\label{appendice}
Consider a Galois extension $L$ of $k$ of group $G$; then $Y= \Spec (L) \to \Spec (k)$ is a right $G$-torsor. The left torsor $Y^0$ is $Y^0=\Spec (L)$ where $L$ is endowed with the right action of $G$ given by $\alpha \star \sigma = \sigma ^{-1} . \alpha $. Then $\Phi ^Y(Y^0)=Y \wedge ^G Y^0= \Spec (L \otimes _k L)^G$ where this denotes the fixed elements of $L \otimes _k L$ under the diagonal action of $G$.

Let us explicit this diagonal action. The \'etale $k$-algebra $L \otimes _k L$ is isomorphic to $L^{(G)}$ on the following way: 

$$\alpha \otimes _k \beta \to (\alpha \sigma ( \beta ) ) _{\sigma \in G}.$$
Let $ \tau \in G$; the image of $\tau (\alpha ) \otimes _k \tau (\beta )$ under this isomorphism is $(\tau (\alpha ) \sigma (\tau ( \beta )) _{\sigma \in G}= ( \tau (\alpha (\tau ^{-1} \sigma \tau ) (\beta )) )) _{\sigma \in G}$. One deduces that the image of an element $(x_\sigma ) _{\sigma \in G}$ by the diagonal action of $\tau \in G$ is $(y_\sigma ) _{\sigma \in G}$ where $$ y_\sigma = \tau (x_{\tau ^{-1} \sigma \tau}).$$ 

The fixed elements under the action of $G$ are the tuplets $(x_\sigma ) _{\sigma \in G}$ satisfying
$$\forall \sigma \quad  \forall \tau \quad x_\sigma = \tau (x_{\tau ^{-1} \sigma \tau}).$$

Consider the partition  $G = \cup _{1\leq i \leq d}C_i$ in conjugation classes and for every $i$, a representative $\sigma _i$ of $C_i$ and $Z_i=Z(\sigma _i)$ the centralizer of $\sigma _i$. From the above condition one deduces that $x_{\sigma _i} \in L^{Z_i}$ (where $L^{Z_i}$ denotes the field of fixed elements under $Z_i$). For any other element $\sigma \in C_i$ of the form $\sigma = \theta \sigma _i \theta ^{-1}$ ($\theta $ is defined up to right multiplication by an element of $Z_i$) we have $x_\sigma =\theta (x_{\theta ^{-1} \sigma \theta } ) = \theta (x_{\sigma _i }) \in L ^{Z(\sigma )}=\theta (L^{Z_i})$.

From that one deduces the inclusion $(L\otimes _k L ) ^G \subset \prod _{1 \leq i \leq d}L^{Z_i}$ and, considering the dimensions over $k$, $$(L\otimes _k L ) ^G \simeq \prod _{1 \leq i \leq d}L^{Z_i}.$$
So sections over $k$ of $\underline {\Isom }_G (L, L)$ are in one to one correspondence with elements of $Z(G)$. If for instance $G$ is abelian, all sections of $\underline {\Isom }_G (L, L)$ are defined over $k$ and $\underline {\Isom }_G (L, L)(k)= G$.

Finally $\Phi ^Y(Y^0)=Y \wedge ^G Y^0= \Spec (L \otimes _k L)^G$ is a right torsor under the inner form $H$ of $G$ defined by the torsor $Y$. Suppose that this torsor $Y$ is defined by a morphism $\Phi : \Gal ( \bar k /k) \to G$, $H$ can be described by the action of $\Gal (\bar k /k)$ on its $\bar k$-points i.e. the elements of the abstract group $G$ by the formula
$$\forall \sigma \in \Gal ( \bar k /k) \quad \forall g \in G \quad \quad \sigma \star g = \Phi ( \sigma ) g \Phi (\sigma )^{-1}.$$

\section{Arithmetic torsion}
\medskip

Let us give an affine $S$-group $G$. We will call {\it an arithmetic torsor} on $f : X\to S$ a torsor coming  by pull back by $f$ from a torsor over $S$. In this section we will consider the particular case of the twisting of an arithmetic $G$-torsor by a $G$-torsor $Q \to X$ over $X$. In other words $\Isom (P_X,Q)= \Phi ^Q (P_X^0)=Q \wedge ^G P_X^0$, where $P \to S$ is a $G$-torsor over $S$. In this situation, $P^0$ being locally trivial on $S$, $\Phi ^Q (P_X^0)$ is locally on $S$ isomorphic to $Q$. We will say that it is a model of $Q$.  

The following statements are immediate consequences of the properties of the twisting operation. 

\begin{prop}\label{specialisation} Let $u: U \to X$ be a $S$-morphism and $s:U \to S$ the composed morphism $s=f\circ u$. The following conditions are equivalent :
\begin{enumerate}
\item $u^\ast Q \simeq s^\ast P$ ;
\item there exists a section $U \to u^\ast (\Phi ^Q (P_X^0))$.

In the case $U= \Spec (k)$, where $k$ is a ring, the preceding conditions are equivalent to the following
\item  $\Phi ^Q (P_X^0)(k) \not= \emptyset$.
\end{enumerate}

\end{prop}
 \proof It is an immediate consequence of Corollary \ref{isomorphisme} and of base change properties. \endproof

\medskip

One can state consequences of Proposition \ref{specialisation} in particular situations. 

\begin{prop}\label{henselian} Let $S= \Spec (R)$ where $R$ is an henselian ring whose we will denote the generic and special points $\eta $ and $s$, $k$ the field of fractions of $R$ and $\kappa $ its residue field. Let $X \to S$ be a smooth $S$-scheme, $G \to S$ an \'etale group scheme, $PÊ\to S$ a right $G$-torsor, $Q \to X$ a right $G_X$-torsor. One assumes there exists a special point $x : \Spec (\kappa ) \to X$ and a $\kappa $-rational point $y:\Spec (\kappa ) \to \Phi ^Q (P^0_X)$ over $x$. 

Then there exists generalisations $\tilde x : S \to X$ of $x$  and these $\tilde x$ satisfy $\tilde x ^\ast Q \simeq P$. In particular $\tilde x _\eta ^\ast Q \simeq P_\eta$.

\end{prop}

\proof The existence of generalizations $\tilde x : S \to X$ of $x$ is due to the fact that $X \to S$ is smooth using Hensel Lemma. On the other hand $\Phi ^Q (P^0_X)$ is \'etale over $X$ as $P_X \to X$ is \'etale and the section $\tilde x$ lifts to a section $\tilde y : S \to \Phi ^Q (P^0_X)$ : indeed the following diagram is cartesian $$\xymatrix{x^\ast \Phi ^Q(P_x^0) \ar[d] \ar[r] &\tilde x ^\ast \Phi ^Q(P_X ^0)\ar[d]\\
\Spec (\kappa ) \ar[r]_s & S\\
}$$ 
where the vertical maps are \'etale; the sections of the left vertical map which correspond to $\kappa$-points of $\Phi ^Q(P_X^0)$ over $x$ are in one to one correspondence with sections of the right vertical map which correspond to $S$-points of $\Phi ^Q(P_X^0)$ over $\tilde x$ (\cite{Milne}, Chapter 1, section 4, Proposition 4.4). Finally according to Proposition \ref{specialisation} the $S$-points of $\Phi ^Q(P_X^0)$ over $\tilde x$ correspond to isomorphisms $\tilde x ^\ast Q \simeq P$. By restriction to the generic fiber one deduces an isomorphism $\tilde x_\eta ^\ast Q \simeq P_\eta $. \endproof

\begin{prop}\label{henselian-2} Let $S= \Spec (R)$ where $R$ is an henselian ring whose we will denote the generic and special points $\eta $ and $s$, $k$ the field of fractions of $R$ and $\kappa $ its residue field. Let $X \to S$ be a smooth $S$-scheme, $G \to S$ an affine smooth group scheme, $PÊ\to S$ a right $G$-torsor, $Q \to X$ a right $G_X$-torsor. One assumes there exists a $\kappa $-rational special point $y:\Spec (\kappa ) \to \Phi ^Q (P^0_X)$.

Then there exists a $S$-point $\tilde x : S \to X$ such that $\tilde x ^\ast Q \simeq P$. In particular $\tilde x _\eta ^\ast Q \simeq P_\eta$.
\end{prop}
\proof The proof is similar to that of preceding Proposition. The scheme $\Phi ^Q (P_X^0)$ is smooth over $X$ and thus over $S$ and Hensel's Lemma insures the existence of a section $\tilde y : S \to \Phi ^Q (P_X^0)$ specializing at $y$. If $\tilde x : S \to X$ is its image in $X$, according to Proposition \ref{specialisation}, $\tilde x ^\ast Q \simeq P$.
\endproof

When the residue field $\kappa$ of $R$ is finite, one may apply Lang Weil estimates for the number of $\kappa$-rational points of a $\kappa $-variety (\cite{Lang-Weil}) and apply Propositions \ref{henselian} and \ref{henselian-2} to insure the existence of $S$-points $\tilde x : S \to X$ such that $\tilde x ^\ast Q \simeq P$.  In order to do this we will have to check whether the special fiber of $\Phi ^Q (P_X^0)$ has $\kappa $-rational geometrically irreducible components, a question that we will address in section \ref{section-specialization}.

\section{Twisting a ramified cover}
A ramified cover $F:Y \to X$ is a finite faithfully flat morphism whose restriction to some dense open $U\subset X$ is \'etale. Let $S$ be a scheme, $X,Y$ be $S$-schemes, $G \to S$ an \'etale goup scheme and $F: Y\to X$ be a $S$-morphism which is a ramified cover endowed with a compatible action of $G_X$ on $Y$ which makes the restriction of $F$ to $U$ a left $G$-torsor. Let $P\to S$ be a right $G$-torsor. According to Theorem \ref{thgeneral} the twisted object $\Phi ^{P_X}(Y)$ is well defined. The restriction of $\Phi ^{P_X}(Y)$ to $U$ is isomophic to $\Phi ^{P_U} (Q)$, where $Q\to U$ is the left $G$-torsor restriction to $U$ of $F: YÊ\to X$ and is in particular \'etale. As $P$ is locally trivial for the \'etale topology, $\tilde F: \tilde Y=\Phi ^{P_X}(Y)\to X$ is \'etale locally on $S$ isomorphic to $Y \to X$. In particular $\Phi ^{P_X}(Y)\to X$ is finite flat ( \cite{EGAIV-2} Proposition 2.7.1). So $\tilde F: \Phi ^{P_X}(Y)\to X$ is a ramified cover. As a consequence of Corollary \ref{specialisation} one gets the following result.

\begin{prop}\label{revetement-corps} There exists a ramified cover $\tilde F :\tilde Y \to X$ \'etale locally isomorphic on $S$ to $F : Y \to X$ such that for any $S$-scheme $t :TÊ\to S$ there exists a $T$-point $y: T \to \tilde Y$ over a point $x \in U(T)$ if and only if $x^\ast Y \simeq t^\ast P^0$ as $G_T$-left torsors.
\end{prop}

In the particular case where $S= \Spec (k)$ of a field $k$ one obtains the following consequence :

\begin{cor} 
Let $F: YÊ\to X$ be a ramified cover over $k$ endowed with a compatible action of $G_X$ on $Y$ which makes the restriction of $F$ to the complement of the branch locus of $f$ a left $G$-torsor. Let $P\to \Spec (k)$ be a right $G$-torsor.  Then there exists a model $\tilde F : \tilde Y \to X$ over $k$ de $F_{\bar k} : Y_ {\bar k} \to X_{\bar k}$ satisfying the following property : for any extension $k'$ of $k$ and for all unramified $x \in X(k')$ the fiber of $F$ at $x$ is isomorphic to the $G$-torsor $P^0_{k'}$ if and only if the fibre at $x$ of $\tilde F$ has a $k'$-rational point. This statement applies in particular to Galois ramified covers.
\end{cor}

\proof It is a consequence of Corollary \ref{specialisation} applied to the $G$-torsor $Q$ obtained from $F: Y \to X$ by removing the branch locus from $X$. The model $\tilde F : \tilde Y \to X$  is the unique finite cover of $X$ whose restriction to the complement of the branch locus of $F$ is isomorphic to $\Phi ^{P_X} (Q)$ (recall that $\Phi ^{P_X} (Q)\simeq (\Phi ^{Q^0} (P_X^0))^0$). \endproof

\begin{prop} \label{revetement-anneau} Let $S= \Spec (R)$ where $R$ is a discrete valuation ring, $\eta $ and $s$ the generic and special points, $k$ the field of fractions of $R$ and $\kappa $ its residue field. Let $X \to S$ be a proper $S$-scheme which we assume to be normal and connected, $G$ be a constant finite group, $PÊ\to S$ be a right-$G$-torsor for the \'etale topology, $Z\to X_\eta $ be a Galois ramified cover of group $G$, with $Z$ normal. We assume that the normalization $F:Y \to X$ of $X$ in $Z\to X_\eta $ is \'etale outside a closed $S$-subscheme $D \not= X$. There exists a model $\tilde F: \tilde Y \to X$ (in the \'etale local sense over $S$), such that

\begin{enumerate}
\item if there exists a $k$-rational point $ y \in \tilde Y_\eta (k)$ over an unramified point $x_\eta \in (X\setminus D)( k)$ the fiber at $x_\eta $ of $Z \to X_\eta $ is isomorphic to the $G$-torsor $P^0 _\eta $. 

\item if moreover the unique extension $x \in X( R)$ of $x_\eta $ doesn't meet $D$, $x^\ast Y\simeq P^0$.

\item one supposes here that $R$ henselian and that $X$ is smooth over $S$ ; if $\tilde Y _s$ has a $\kappa $-rational point $v \in  \tilde Y_s (\kappa )$ over a point $u \in X_s (\kappa )$ not belonging to $D$, there exists a section $ x \in (X\setminus D)(R)$ extending $u$ such that $x^\ast Y\simeq P^0$.
\end{enumerate}

\end{prop}

\proof As a consequence of the hypotheses the generic fiber of $Y \to X$ is isomorphic to $Z \to X_\eta $. The action of $G$ on $Z$ extends in an action of $G$ on $Y$  due to the normality of $Y$. Denote $U= X\setminus D$, $Q = F^{-1}  (U)$,  by hypothesis $Q \to U$ is \'etale. As $U$ is normal, the same is true for $Q$ (\cite{EGAIV-2}, Corollaire 6.5.4). The same arguments show that $Q \times _{X\setminus D} Q$ and $Q\times _S G$ are normal. From that it follows that the isomorphism $f:Q_\eta \times _k G \simeq Q_\eta \times _{X_\eta \setminus D_\eta}Q_\eta $ and its inverse $g=f^{-1}$ extend in morphisms $\bar f: QÊ\times _S G \to QÊ\times _{X\setminus D} Q$ and $\bar g :  QÊ\times _{X\setminus D} Q \to QÊ\times _S G$ over $U$ and the restrictions of $\bar g\circ \bar f$ and $\bar f \circ \bar g$ to the generic fibers are the identity. Then $\bar f$ and $\bar g$ are isomorphisms inverses  from eachother and $Q \to U$ is a $G$-torsor.

One can consider the twisted object $\Phi ^{P_U} (Q)\to U$ which is locally isomorphic to $Q\to U$ for the \'etale topology and thus \'etale. Moreover $\Phi ^{P_U} (Q)\to U$ is the restriction to $U$ of $\tilde Y=\Phi ^{P_X} (Y)\to X$.
\begin{enumerate}

\item the point $ y \in \tilde Y_\eta (k)$ belongs to $\Phi ^{P_U} (Q)(k)$ over $x_\eta \in U(k)$. The conclusion follows from \ref{specialisation}.

\item the restriction to $U$ of $\tilde Y \to X$ is finite \'etale and the unique section $\tilde x \in U(R)$ lifts to a section $\tilde y \in \tilde Y (R)$ and belongs in fact to $\Phi ^{P_U} (Q)(R)$. The conclusion follows from Proposition \ref{specialisation}. 

\item As $U \to S$ and $\Phi ^{P_U} (Q)\to S$ are smooth the point $v \in  \tilde Y_s (\kappa )$ extends in a section $y \in \Phi ^{P_U} (Q)(R)$ over a section $ x \in U(R)$ which satisfy according to Proposition \ref{specialisation} $x^\ast Q \simeq x^\ast Y\simeq P^0$.
\end{enumerate}

\endproof
An example of situation where Proposition \ref{revetement-anneau} apply is given by the following statement.

\begin{cor} \label{cor-revetement-anneau} Let $S= \Spec (R)$ where $R$ complete discrete valuation ring, $\eta $ and $s$ the generic and special points, $k$ the field of fractions of $R$ and $\kappa $ its residue field. Let $X \to S$ be a smooth proper relative curve over $S$, $G$ a finite constant group, $PÊ\to S$ a right $G$-torsor for the \'etale topology, $Z\to X_\eta $ a Galois ramified cover of group $G$. One assumes that the normalization $F:Y \to X$ of $X$ in $Z\to X_\eta $ has no vertical ramification. Then conclusions of Proposition \ref{revetement-anneau} hold. 
\end{cor}

\proof Under these hypotheses  the morphism $F : Y \to X$ is flat and defines a ramified cover along a divisor $D$ whose components are the closure in $X$ of the branch points of the cover $Z\to X_\eta $. One can apply Proposition \ref{revetement-anneau}.
\endproof

\begin{rem} When the center $Z(G)$ of $G$ is trivial, it follows from \cite{Beckmann}, Propsition 2.3 that if the residue characteristic doesn't divide the order of $G$ and the distinct branch points don't meet on the special fiber, the cover $F: Y \to X$ has no vertical ramification. One may apply \ref{cor-revetement-anneau} in this situation. 
\end{rem}

Let us end this section by a statement which illustrates how Lang-Weil's estimates (in this instance Riemann hypothesis in function fields) can be used in this context.

\begin{prop}\label{Lang-Weil-1} Let $R$ be a henselian discrete valuation ring with finite residue field $\kappa $ and with fraction field $k$, $X \to S= \Spec (R)$ be a smooth proper $R$-curve with $\h ^0 (X , \dO _X)=R$, $f: Y \to X$ a ramified cover with no vertical ramification (finite \'etale over some open $U\subset X$ which surjects onto $S$), $G \to S$ an \'etale finite group scheme acting on $Y$ over $X$ such that the restriction $Q\to U$ of $f$ to $U$ is a left $G$-torsor.  Assume $Y \to S$ to be smooth and $\h ^0 (Y, \dO _Y)=R$.

Then there exists a constant $C>0$ depending on the degree $d$ of the cover $f$, the genus $g$ of the fibers of $Y \to S$ and the number $r$ of the branch points, such that for any finite integral ring extension $R \subset R'$ of residue field $\kappa '$, with $[\kappa ' : \kappa ] \geq C$ and any right $G$-torsor $P \to \Spec (R')$, there exists unramified $R'$-points $x: \Spec (R') \to X$ such that $x^\ast Y \simeq P^0$ as left $G$-torsors. 
\end{prop}

\proof It follows from the hypothesis that $Y\to S$ has smooth and geometrically connected fibers and thus geometrically irreducible fibers. As $\tilde Y = \Phi ^{P_X} (Y)\to S$ is \'etale locally isomorphic to $Y\to S$, $\tilde Y \to S$ has geometrically irreducible fibers. In particular the special fiber $\tilde Y _s \to X_s$ is a smooth ramified cover with less than $rd$ ramification points and $\tilde Y_s$ is geometrically irreducible. For any finite extension of the residue field $\kappa \subset \kappa '$, the number $N_{\kappa '}$ of $\kappa '$-points on $\tilde Y_s$ satisfies the inequality

$$| N_{\kappa '} -(|\kappa ' |+1)| \leq 2g \sqrt {|\kappa ' |}.$$ So for $|\kappa ' |$ large enough (depending on $g,r,d$) $N_{Ê\kappa '} > rd$, and there exists unramified points $\kappa '$-points $v$ on $\tilde Y_s$. Let $U$ be the complement of the branch locus and $Q = f^{-1} (U)$; then $Q \to U$ is a left $G$-torsor and $v$ extends in a $R'$-point $y: \Spec (R') \to \Phi ^{P_U}Q$ for any discrete valuation ring extension $R'$ of residue field $\kappa '$. Let $x : \Spec (R') \to U \subset X$ the image of $y$. As in point (3) of Proposition \ref{revetement-anneau}, $x ^\ast Y \simeq P^0$. \endproof

\section{Specialization of a cover}\label{section-specialization}
Let $k$ a field and $X \to \Spec (k)$ a proper $k$-sheme. Let $F : Y \to X$ be a ramified cover and $U \subset X$ be a dense open subscheme such that the restriction $Q=F^{-1} (U) \to U$ of $F$ above $U$ is finite \'etale. We assume $X$ to be geometrically normal and geometrically connected and $Y$ geometrically normal (geometrically normal resumes to normal if the base field is perfect, see \cite{ST}, Lemma 10.151.1 ). The open subscheme $U$ is obviously geometrically normal. On the other hand, as for any field extension $k \subset k'$, $X_{k'}$ is the normalization of itself in $U_{k'} \hookrightarrow X_{k'}$, $U$ is also geometrically connected. As $F$ is faithfully flat, it is open, and $Q$ is a dense open subscheme of $Y$. This is true for any base field extension, and thus if $Q$ is geometrically connected, $Y$ is geometrically connected. On the other hand, for any field extension $k \subset k'$, $Y_{k'}$ is the normalization of $X_{k'}$ in $Q_{k'} \to X_{k'}$, and thus, if $Y$ is geometrically connected, $Q$ is geometrically connected. 

Describing finite \'etale covers of $U$ in terms of morphisms $\pi _1 (U, \bar x) \to \dS _d$, one sees that $Q$ is geometrically connected if and only if there isn't a non trivial finite field extension $k \hookrightarrow L$ such that $Q \to U$ factors through $Q \to U_L \to U$. In particular if $Q$ is not geometrically connectedl $Q(k)=\emptyset $

In the case $F : Y \to X$ is Galois of group $G$, then $F_U:Q=F^{-1} (U) \to U$ is a left $G$-torsor under the constant group $G$. The above remarks can easily be formulated in terms of morphisms of \'etale fundamental groups. The \'etale cover $Q \to U$ is described by a surjective morphism $ \Phi : \pi _1 (U, \bar x) \to G$ (where $\bar x$ refers to a geometric point of $U$). This morphism $\Phi $ inserts in the following commutative diagram where the vertical maps are surjective

\begin{equation}\label{2}
\xymatrix{1 \ar[r] & \pi _1 (U_{\bar k} , \bar x) \ar[r] \ar[d]^\varphi & \pi _1 (U, \bar x) \ar[r]^u\ar[d]^{\Phi } & \Gal ( \bar k /k) \ar[r] \ar[d] & 1\\
1\ar[r] & \bar G \ar[r] & G \ar[r]^v& \Gal ( L /k) \ar[r] & 1\\}
\end{equation}
where $L$ is the scalar extension in the covering. This means that the \'etale cover $Q \to U$ factors through $Q \to U_L \to U$. The condition $L =k$ (or equivalently $\bar G=G$) is equivalent to the condition that $Q$ is geometrically connected or equivalently that $Y$ is geometrically connected. We conclude that a necessary condition for $Q$ to have a $k$-rational point is that $Q$ is geometrically connected.

\medskip

Coming back to the problem to know, given a left $G$-torsor $P= \Spec (K) \to \Spec (k)$ (where $K$ is a finite \'etale $k$-algebra), whether the specialization of $Q \to X$ at some rational point $x \in U(k)$ is isomorphic to $P\to \Spec (k)$, Proposition \ref{revetement-corps} gives an answer in terms of $k$-rational unramified point in $\tilde Y$. In what follows we won't assume $Y$ to be geometrically connected. To apply Proposition \ref{revetement-corps} one has to be able to describe the connected components of $\tilde Y$ and to check whether they are geometrically connected. This the aim of this more technical section.
\medskip
In what follows we will assume that the following condition is satisfied :

$(\star )$ \hskip0,5cm {\it the quotient of the left $G$-torsor $P= \Spec (K) \to \Spec (k)$ by $\bar G$ is isomorphic to the $G / \bar G$-left torsor $\Spec (L) \to \Spec (k)$}

\medskip

It is an obvious necessary condition for $Q\to X$ to having a rational specialization isomorphic to $P \to \Spec (k)$. Suppose that $P\to \Spec (k)$ is described by the morphism $\Psi : \Gal (\bar k /k) \to G$ whose kernel $\Gal (\bar k /N)$ defines a Galois extension $N$ of $k$. If  $H=\Gal (N/k)$ one may suppose $L\subset N$, $L=N^{\bar G \cap H}$ and $\bar G H=G$.

We have seen that a necessary and sufficient condition for the problem to have an affirmative answer is that $Q^0\wedge ^{G} P_U$ has a $k$-rational point $y$; if $x \in X(k)$ is the image of $y$, the fiber of $Q$ at $x$ will be isomorphic to $P$. We are in the situation of Proposition \ref{quotient-partiel}: the $G$-torsors $Q$ and $P_U$ have a common quotient $R_U$. The twisted object $R^0\wedge ^{G/\bar G} R=\underline {\Isom }_{G/\bar G} (R,R)$ is a trivial torsor whose sections over $k$ are in one to one correspondence with $Z(G/\bar G)$. It follows that if $X$ has $k$-rational points (otherwise the question is empty) the sections $X \to \underline {\Isom } _{G/\bar G} (R_X,R_X)$ are themselves in one to one correspondence with $Z(G/\bar G)$. Let $s$ be such a section. Let us recall the diagram
\begin{equation}
\xymatrix{
s^\ast \underline {\Isom } _{G} (Q,P_U) \ar[r] \ar[d] &\underline {\Isom } _{G} (Q,P_U) \ar[d] \\
 X \ar@/_/@{-->}[r]_s& \underline {\Isom } _{G/\bar G} (R_U,R_U)\ar[l] \\
}
\end{equation}

and if one pulls it at the rational point $x \in U(k)$ one gets
\begin{equation}
\xymatrix{
x^\ast (s)^\ast \underline {\Isom } _{G} (x^\ast Q,P) \ar[r] \ar[d] &\underline {\Isom } _{G} (x^\ast Q,P) \ar[d] \\
 \Spec (k) \ar@/_/@{-->}[r]_{x ^\ast (s)} & \underline {\Isom } _{G/\bar G} (R,R)\ar[l] \\
}
\end{equation}
The sections $\Spec (k) \to \underline {\Isom } _{G} (x^\ast Q,P)$ induce sections $x ^\ast (s): \Spec (k) \to \underline {\Isom } _{G/\bar G} (R,R)$ for some $s$ corresponding to some element of $Z(G/\bar G)$ and for a given $s$ are in one to one correspondence with sections $\Spec (k) \to x^\ast (s)^\ast \underline {\Isom } _{G} (x^\ast Q,P)$, in other words $k$-rational points of $s^\ast \underline {\Isom } _{G} (Q,P_U) $ above $x$. We state the result:

\begin{prop}\label{decomposition} Assume that condition $(\star )$ is fulfilled. With the preceding notation let $\{ s_\gamma \} _{\gamma \in Z(G /\bar G)}$ be the set of sections $\Spec (k) \to \underline{\Isom }_G (R,R)$. Then for all $\gamma \in  Z(G /\bar G)$, $s_\gamma ^\ast (Q\wedge ^GP^0_{U})$ is geometrically connected. Moreover there exists a unramified rational point $x \in X(k)$ such that the fiber at $x$ of $Y$ is isomorphic to the $G$-torsor $P$ if and only if there exists $ \gamma \in Z(G/\bar G)$ such that $s_\gamma ^\ast (Q\wedge ^GP^0_{U})(k) \not= \emptyset$.
\end{prop}

\proof The only thing to prove is that the $s_\gamma ^\ast (Q\wedge ^GP^0_{X})$'s are geometrically connected. Recall the cartesian diagram of Proposition \ref{quotient-partiel}. 
\begin{equation}
\xymatrix{
\Isom _{\bar G,X_L}(Q, (1_X \times \gamma)^\ast P_X )\ar[r] \ar[d] &        s^\ast \Isom _{G} (Q,P_X) \ar[r] \ar[d] &\Isom _{G} (Q,P_X) \ar[d] \\
X_L \ar[r] & X \ar@/_/@{-->}[r]_{s_\gamma }& \Isom _{G/\bar G} (X_L,X_L)\ar[l] \\
}
\end{equation}
But $(1_X \times \gamma)^\ast P_X \simeq (\gamma ^\ast P)_X$ where $\gamma ^\ast P \to \Spec (L)$ is a $\bar G$-torsor. And $\Isom _{\bar G,X_L}(Q, (1_X \times \gamma)^\ast P_X )\simeq Q\wedge ^{\bar G} (\gamma ^\ast P)_{X_L}^0$ is isomorphic over $X_{Ê\bar k}$ to $Q\times _L{\bar k}$ which is connected as $L$ is the extension of scalars contained in $Q$. 

\endproof

\begin{rem} Instead of sections over $k$ one can more generally consider sections over a finite extension of $k$. For instance if one works over $L$, sections $\Spec (L) \to \underline{ \Isom} _G (R,R)$ are in one to one correspondence with $G/\bar G$ (the inner form of $G/\bar G$ induced by the torsor $R= \Spec (L) \to \Spec (k)$ split over $L$ and is thus isomorphic to $G/\bar G$). Thus over $L$, $\underline{\Isom} _{G} (Q_L,P_{X_L})=Q_L\wedge ^GP^0_{X_L}$ is a disjoint union of the open closed $s_\gamma ^\ast \underline{\Isom} _{G} (Q_L,P_{X_L})$ with $\gamma $ running in $G/\bar G$. The situation reduces to that of the regular $G$-cover $Q_L \to X_L$.
\end{rem}

\medskip

When $k$ is a {\it pseudo algebraically closed} field (abbreviated by PAC) the conclusion of Proposition \ref{decomposition} always holds. Recall that a field $k$ is a PAC field if every geometrically irreducible variety defined over $k$ has $k$-rational points (see \cite{Fried-Jarden}). As in Proposition \ref{decomposition} the components $s_\gamma ^\ast (Q\wedge ^G P^0_U)$ are geometrically irreducible, one gets the following corollary.

\begin{cor} Assume that condition $(\star )$ is fulfilled. Then if $k$ is a PAC field there are infinitely many unramified rational points $x \in X(k)$ such that the fiber at $x$ of $Q$ is isomorphic to the $G$-torsor $P$.
\end{cor}

\medskip

One can generalize Proposition \ref{decomposition} to the situation of schemes over a discrete valuation ring instead of a field. To avoid confusion in the notation, let us call $A$ this discrete valuation ring; we denote as usually $k$ its field of fractions and $\kappa $ is residue field. Let $X \to \Spec (A)$ be a faithfully flat and proper $A$-scheme that we will assume to be integral and normal. Let $Z \to X_\eta $ be a Galois ramified cover of the generic fiber of Galois group $G$ and geometric Galois group $\bar G <G$. Let us assume that the normalized $f:Y\to X$ of $X$ in $Z$ is flat \footnote{This will be always the case for a regular scheme $X$ of relative dimension $1$ .}. Assume also that $f: Y \to X$ has no vertical ramification (this means that its restriction to the special fiber is a ramified cover). Let $\Spec (L) \to \Spec (k)$ be the extension of scalars in the Galois cover $Z \to X_\eta $ ($L$ is a Galois extension of $k$ of group $G/\bar G$) and $A_L$ the integral closure of $A$ in $L$ that we assume \'etale over $A$. We suppose that there exists an $A_L$-point $x : \Spec (A_L) \to X_{A_L}$. Let $P \to \Spec (A)$ be an \'etale $G$-torsor such that the generic fiber $P_\eta \to \Spec (k)$ factors through $P_\eta \to \Spec (L) \to \Spec (k)$. We are in the situation of Proposition \ref{quotient-partiel}. 

\begin{prop} \label{decomposition-anneau}Under these hypothesis
\begin{enumerate} 
\item  $R=\Spec (A_L) \to \Spec (A)$ is an \'etale $G/\bar G$-torsor and there are factorizations $P \to \Spec (A_L) \to \Spec (A)$ and $Y \to X_{A_L} \to X$; 
\item $\Isom _{G/\bar G} (X_{A_L},X_{A_L}) =\Isom _{G/\bar G} (R,R) = \h^0 (\Spec (A),\underline{ \Isom} _{G/\bar G} (R,R) )\simeq \Isom _{G/\bar G} (R_\eta ,R_\eta )\simeq Z(G/\bar G)$;
\item for all $ \gamma \in Z(G/\bar G)$, if we denote by $s_\gamma \in \Isom _{G/\bar G} (X_{A_L},X_{A_L}) =\Phi ^{R_X^0}(R_X)$ the corresponding element, $s_\gamma ^\ast (\Phi ^{P^0_X} (Y))$ has geometrically connected fibers.
\end{enumerate}
\end{prop}
\proof The first assertion comes from the normality of $P$ and $Y$. In the second assertion the first equality comes from the existence of a point $\Spec (A_L) \to X$. Other equalities are clear. As $L$ is the constant field extension in $f_\eta : Y_\eta \to X_\eta $, one gets $\h^0 (Y_\eta , \dO _{Y_\eta })= L$. On the other hand $\h^0 (Y_\eta , \dO _{Y_\eta }) \simeq \h ^0 (Y, \dO _Y) \otimes _A k$ (cf. \cite {Liu}, section 5, Ex. 1.16, p. 174). As $\h^0 (Y, \dO _Y)$ has no $A$--torsion, so $A_L \subset \h^0 (Y, \dO _Y) \subset \h ^0 (Y, \dO _Y) \otimes _A k=L$ and as $\h^0 (Y, \dO _Y)$ is a finite $A$-algebra by Serre's theorem, it is the integral closure of $A$ in $L$: $\h^0 (Y, \dO _Y)=A_L$. It follows from this fact that $Y \to X_{A_L}$ has geometrically connected fibers. As in the proof of Proposition \ref{decomposition} one uses Proposition \ref{quotient-partiel} to show that $s_\gamma ^\ast (\Phi ^{P^0_X} (Y))\to X$ is \'etale locally over $\Spec (A)$ isomorphic to $Y \to \Spec (A_L)$ and thus has geometrically connected fibers. 
\endproof
\medskip

In the situation of Proposition \ref{decomposition-anneau}, let us give an example of application of the Lang-Weil estimates.

\begin{prop}\label{Lang-Weil-2} Under the hypothesis of Proposition \ref{decomposition-anneau}, assume moreover that $X$ and $Y$ are smooth relative curves over $\Spec (A)$ and that $A$ is henselian with finite residue field $\kappa $. There is a constant $C>0$ depending on the degree of the covering $Y \to X $, the genus of the fibers of $Y $ and the number of branch points, such that for any discrete valuation ring extension $A'$ of $A$ whose residue field $\kappa '$ satisfies $[\kappa ' : \kappa ]\geq C$ there are unramified $A'$-points $x : \Spec (A') \to X$ such that $x^\ast Y \simeq P_{A'}$ as left $G$-torsors.  
\end{prop}

\proof The proof is similar to the proof of Proposition \ref{Lang-Weil-1} taking advantage of the fact that $Y \to X$ is smooth and $A$ is henselian. One uses Lang-Weil estimates \cite{Lang-Weil} to insure the existence of $\kappa '$-unramified points on the special fiber of $s_\gamma ^\ast (\Phi ^{P^0_X} (Y))$ for $[\kappa ' : \kappa]$ large enough. \endproof
\medskip

The viewpoint of \cite{Debes-Legrand} is different : for a Galois extension $N$ of $k$ the authors ask if there are specializations of the \'etale covering $Q$ of $U$ at unramified points $x \in X(k)$ isomorphic to a disjoint union of $(G:H)$ copies de $N$ forgetting the action of $G$. So they have to consider all $G$-torsors over $\Spec (k)$ associated to injective morphisms $H \hookrightarrow G$ whose composition with the canonical surjection $G \to G/\bar G$ is surjective. One is lead to look for $k$-rational points on a family of schemes indexed by embedding $H \hookrightarrow G$ whose composition with the canonical surjection $G \to G/\bar G$ is surjective. Each of these schemes is the twisted object $Q\wedge P^0_X$ where $P$ runs among the above mentioned $G$-torsors. The answer to the question can be formulated as follows.

\begin{prop} Let $N$ be a Galois extension of $k$ of group $N$. Consider the set $J$ of embeddings $j : H \hookrightarrow G$ whose composition with the canonical surjection $G/\bar G$ is surjective, modulo automorphisms of $G$ fixing $H$. For each $j \in J$ denote by $K_j$ the finite \'etale $k$-algebra such that $P_j = \Spec (K_j) \to \Spec (k)$ is the $G$-torsor contacted product of the $H$-torsor $\Spec (N) \to \Spec (k)$ with the embedding $j:H \hookrightarrow G$. The following conditions are equivalent :

\begin{enumerate}
\item There exists $x \in X(k)$ such that the fiber at $x$ of $Q$ is isomorphic to a disjoint union of $(G:H)$-copies of $\Spec (N) \to \Spec (k)$.
\item There exists $j \in J$ such that $(Q\wedge ^G P_{j,X}^0 ) (k) \not= \emptyset$.
\item There exists $j\in J$ and $\gamma \in Z(G/\bar G)$ such that $s_\gamma ^\ast (Q\wedge ^G P_{j,X}^0 ) (k) \not= \emptyset$.
\end{enumerate}
\end{prop}

\medskip

In \cite{Debes-Legrand}, section 3.2 the authors consider more generally the case of non necessarily Galois covers. Let $F:Y \to X$ be a connected ramified cover of degree $n$ \'etale above a connected dense open $U \subset X$, $\Phi : \pi _1 (U, \bar x) \to \bS _n$ the corresponding morphism, $G<\bS _n$ the image of $\Phi $, $Z \to F^{-1} (U) \to U$ the Galois closure of $F^{-1} (U)Ê\to U$, corresponding to the surjective morphism $\phi :  \pi _1 (U, \bar x) \to G$ and $Q\to U$ the associated $\bS _n$-torsor which is the contracted product of $Z$ by $\bS _n$ via the inclusion $G<\bS _n$. Let $L$ be the extension of scalars in $Z \to U$. We are also given an extension of \'etale $k$-algebras $N'$ of $k$ of degree $n$ corresponding to a morphism $\Psi : \Gal (\bar k /k) \to \bS _n$ and let $H$ be the image de $\Psi $ which is the Galois group of the Galois closure $N$ of $N'$ over $k$ (which can be viewed as the compositum in $\bar k$ of the Galois closures of the components of $N'$). We call $P \to \Spec (k)$ the $\bS _n$-torsor associated to the morphism $\Psi $. The $\bS _n$-torsor $Q$ is the contracted product $Z \wedge ^G \bS _n$ and splits in $(\bS _n : G)$ connected components isomorphic to $Z$. In the same way the $\bS _n$-torsor $P$ is the contracted product $\Spec (N) \wedge ^H \bS _n$ and has $(\bS _n : H)$ connected components isomorphic to $\Spec (N)$. 

The $\bS _n$-torsor $Q\to U$  restriction to $U$ of $Y\to X$ has the following description, $x$ being a $k$-rational point of $U$ and $s_x$ the corresponding section, defined up to conjugation by an element of  $\pi _1 (U_{\bar k}, \bar x)$   :

\begin{equation}
\xymatrix{1\ar[r] &\pi _1 ( U_{\bar k}, \bar x) \ar[r] \ar[d]^{\varphi } & \pi _1 (U, \bar x) \ar[r]_u \ar[d]^\phi & \Gal (\bar k /k) \ar@/_/@{-->}[l]_{s_x}
 \ar[r] \ar[d] ^\Lambda &1\\
1 \ar[r]& \bar G \ar[r] &  G \ar[r]_v \ar[d]^{\nu } & G/\bar G \ar[r] & 1\\
&&\bS _n&&\\
}
\end{equation}

where $v \circ \phi \circ s_x= \Lambda$. The fiber of $Q$ at $x$ is described by
\begin{equation}
\xymatrix{\Gal ( \bar k /k) \ar[r]^{s_x}& \pi _1 (U, \bar x) \ar[r]^\phi &G \ar[r]^\nu &\bS _n\\
}
\end{equation}

We write $\Phi = \nu \circ \phi$ and $\Phi _x = \nu \circ \phi \circ s_x$. On the other hand the $\bS _n$-torsor $P$ has the following description :

\begin{equation}
\xymatrix{\Gal ( \bar k /k) \ar[r]^\psi &  H \ar[r]^\mu &\bS _n\\}
\end{equation}
Let us denote $\Psi = \mu \circ \psi$. The specialization of $Y\to X$ at $x$ is isomorphic to $\Spec (N') \to \Spec (k)$ if and only if up to conjugation by an element of $\bS _n$, $\nu \circ \phi \circ s_x= \mu \circ \psi $. Identifying $H$ and $G$ with their images in $S_n$ this implies the existence of an embedding $\eta : H \hookrightarrow G$ such that

$$\nu \circ \eta = \mu \quad {\rm and} \quad \eta \circ \psi = \phi \circ s_x.$$

Conversely suppose that there is an embedding $\eta : H \hookrightarrow G$ such that $\nu \circ \eta = \mu$ and such that the $G$-torsor $P' =\Spec (N) \wedge ^H G\to \Spec (k)$ through the embedding $\eta : H \hookrightarrow G$ is isomorphic to the $G$-torsor $x^\ast Z$ described by the morphism $\phi \circ s_x : \Gal (\bar k /k) \to G$. There exists an element $ \omega \in G$ such that 

$$\forall \gamma \in \Gal ( \bar k /k) \quad \omega (\phi \circ s_x ( \gamma ))\omega ^{-1}= \eta \circ \psi ( \gamma )$$ which implies 

$$\forall \gamma \in \Gal ( \bar k /k) \quad \nu (\omega ) (\nu \circ \phi \circ s_x ( \gamma ))\nu (\omega )^{-1}=\nu \circ  \eta \circ \psi ( \gamma )= \mu \circ \psi ( \gamma )$$ and thus up to conjugation by an element of $\bS _n$, $\nu \circ \phi \circ s_x = \mu \circ \psi $ which means that the specialization of $Y \to X$ at $x$ is isomorphic as cover to $N'$.

One can consider the twisted object $Z\wedge ^G {P'_U}^0$ and the existence of $k$-rationals points on this $k$-scheme over a point $x \in X(k)$ 

is equivalent to the existence of an isomorphism of $G$-torsors $x^\ast Z \simeq P'$. One gets the following statement which is a reformulation of the "Twisting Lemma" 3.4 of \cite{Debes-Legrand}:

\begin{prop} With the preceding notation the finite \'etale $k$-schemes $x^\ast Y \to \Spec (k)$ and $\Spec (N') \to \Spec (k)$ are isomorphic if and only if there exists an embedding $\eta : H \hookrightarrow G$ such that

\begin{enumerate}
\item $\nu \circ \eta = \mu $ ; 
\item the fiber at $x$ of $Z\wedge ^G {P'_U}^0$ has a $k$-rational point where $P' =\Spec (N) \wedge ^H G\to \Spec (k)$ through the embedding $\eta : H \hookrightarrow G$. 
\end{enumerate}
\end{prop}

\end{document}